\documentclass[journal,final]{new-aiaa}
\usepackage{color}
\usepackage{algorithm}
\usepackage{graphicx}
\usepackage{amsmath}
\usepackage{longtable,tabularx}

\newcommand{\R}{\mathbf{R}}
\newcommand{\U}{\mathbf{U}}

\newcommand{\w}{\mathbf{w}}

\newcommand{\e}{\mathbf{e}}

\newcommand{\vv}{\mathbf{v}}

\newcommand{\df}[2]{\dfrac{\partial#1}{\partial#2}}

\graphicspath{{./pic/}}

\title{Simplified Eigenvalue Analysis for Turbomachinery
	Aerodynamics with Cyclic Symmetry}
\author[1,2]{Shenren Xu	\footnote{Associate Professor, School of Power and Energy; shenren\_xu@nwpu.edu.cn}}
\affil[1]{Northwestern Polytechnical University, 127 Youyixi Road, Xi’an 710072, China}
\affil[2]{Yangtze River Delta Research Institute of NPU, Northwestern~Polytechnical~University, Taicang~215400,~P.R.~China}

\begin{document}
\maketitle

\begin{abstract}
Eigenvalue analysis is widely used for linear instability analysis in both
external and internal aerodynamics. It typically involves finding the
steady state, linearizing around it to obtain the Jacobian, and then
solving for its eigenvalues and eigenvectors. When the flow is
modelled with Reynolds-averaged Navier--Stokes equations with a
large boundary-layer-resolving mesh, the resulting eigenvalue
problem can be of very high dimensions,
and is thus computationally very challenging. To reduce the
computational cost, a simplified approach is proposed to compute
the eigenvalues and eigenvectors, by exploiting the cyclic symmetric
nature of annular fluid domain for typical compressors. It is shown
that via a rotational transformation, the Jacobian can be reduced to a
block circulant matrix, whose eigenvalues and eigenvectors then can
be computed using only one sector of the entire domain. This simplified
approach significantly lowers the memory overhead and the CPU time
of the eigenvalue analysis without compromising on the accuracy.
The proposed method is applied to the eigenvalue analysis of an
annular compressor cascade with 22 repeated sectors and it is shown
the spectrum of the whole annulus can be obtained by using the
information of 1 sector only, demonstrating the effectiveness of
the proposed method.
\end{abstract}

\section*{Nomenclature}
{\renewcommand\arraystretch{1.0}
\noindent\begin{longtable*}{@{}l @{\quad=\quad} l@{}}
$A$&Jacobian matrix for the whole annulus\\
$B$   & circulant matrix \\
$b_i$     & entries of a circulant matrix\\
$\e,\hat \e$& basis vector of the original and rotated references\\
$\R$&residual vector\\
$T$&rotational transformation matrix\\
$\U$& flow field vector\\
$\vv,\w$     & eigenvectors \\
$\rho_m$     & $M$th complex roots of 1\\
$\theta$& pitch angle between two neighboring sectors\\
$\lambda$     & eigenvalue\\

\end{longtable*}}

\section{Introduction}
Eigenvalue analysis based on the Reynolds-averaged Navier--Stokes
(RANS) equations is a widely used approach to study the linear instability
of external flow problems~\cite{theofilis2011global} such as the
transonic shock buffet for
airfoil~\cite{Crouch2007Predicting,Crouch2009Origin} and for
wings~\cite{2018arXivTimme,paladini2019transonic,crouch2019global}.
In those works, the RANS equations are linearized about the steady state
and an eigenvalue problem is solved for the resulting large sparse linear system
of equations. The stability of the system is then determined by the existence
of eigenvalues with positive real parts. This approach is shown to be capable
of determining the stability boundary that is consistent with unsteady RANS
approach, but at orders-of-magnitude lower a cost.

Eigenvalue-based linear instability analysis has also been used for
turbomachinery applications, mainly regarding the flow-instability
inception towards the stall boundary of compressors. A general
approach to predict the onset of such instability using the eigenvalue
method is proposed in~\cite{Sun2013A} and has been successfully
applied to both axial and centrifugal compressors~\cite{liu2014basic,sun2016flow}.
Different from the application on airfoil and wing aerodynamics
where the flow is modelled with three-dimensional RANS equations,
these works on compressor stability study use a simplified approach
by modelling the effect of the blades with body forces
and ignoring the circumferential variation of the flow field.
With such simplification, a small eigenvalue problem can be
formulated and solved in the meridional plane to predict the linear
instability of the system. Although the method is shown to predict
the turbomachinery flow-instability onset with plausible accuracy,
ignoring the three-dimensional details naturally brings error that is
not known a priori. Besides, modern compressors are designed with
increasingly higher loading and exhibit strong three-dimensional effect,
which needs to be accounted for when studying the flow instability~\cite{pullan2015origins}.
Furthermore, three-dimensional RANS calculation has become routine practice
for turbomachinery industry and consequently the flow-instability inception
prediction needs to be based on RANS equations as well in order to retain
the same modelling capability.

One key factor limiting the application of the RANS-based eigenvalue analysis
to turbomachinery is its high computational cost. Although the steady state
performance of compressors can be computed using a single passage mesh,
which typically has grid points on the order of one million, whole annulus domain
needs to be employed in order to capture the instability onset. Consequently, a
computational mesh one to three orders of magnitude larger than typical
airfoil/wing applications is needed, and therefore computing eigenvalue
(even a small subset of it) based on RANS equations for compressors
is computationally much more demanding.

Turbomachinery components usually assume cyclic symmetric shape and so is
the flow field inside, and the cyclic symmetry can be exploited in order to simplify
the computation. It is in fact established practice for computing the natural
frequencies and vibrational mode shapes of periodic structures with cyclic
structure~\cite{thomas1979dynamics}.
One first meshes one sector of the entire structure and specifies
the nodal diameter and the natural frequencies and mode shapes for the
specifies nodal diameters can be computed. Looping over all possible
nodal diameters allows one to obtain all natural frequencies and mode shapes
for the entire structure. Throughout the whole process, all computations
are performed based on one sector of the domain, thus requires much
less computational resource and much shorter CPU time, compared to
the whole-annulus computation.
Such simplification approach was revisited in~\cite{Schmid2017Stability}
for flow problems with translational periodicity and was successfully demonstrated
on linear cascade type of problems. However, the methodology proposed
therein is not directly applicable to cases with rotational periodicity, which
are representative of real compressors. Therefore, it is desirable to extend
the methodology to flow problems with rotational periodicity in order to
significantly simplify the eigenvalue computation in order to study
the whole annulus compressor flow instability.

In this work, an extension of the method proposed in~\cite{Schmid2017Stability}
is developed so that it can be directly applied to the eigenvalue calculation for
cyclic symmetric domains using one sector of the geometry only. In the
remaining part of the paper, the theoretical basis for the simplification of the
eigenvalue computation is first reviewed in Sec.~\ref{sec:theory} with
particular focus on the linearized RANS equations. The application of
the proposed method is then applied to the eigenvalue calculation for
an annular compressor cascade case to demonstrate its effectiveness
in reducing the computational cost and results are discussed in Sec.~\ref{sec:result}. Conclusions are given in Sec.~\ref{sec:conclusion}.

\section{Theoretical framework}
\label{sec:theory}
\subsection{Eigenvalue analysis of block circulant matrices}
A circulant matrix $B$ has the form as follows
\begin{equation*}
	B=\left [
	\begin{array}{c c c c }
		b_0       & b_1     &  \cdots &  b_{M-1}\\
		b_{M-1} & b_0       &  \cdots &  b_{M-2}\\
		\vdots   & \vdots & \ddots & \vdots\\
		b_1        & b_2        &  \cdots & b_{0}
	\end{array}
	\right].
\end{equation*}
Due to its cyclic symmetric nature, it can be verified
that matrix $B$ has the following eigenvalues and eigenvectors
\begin{equation*}
\lambda_m= \sum_{k=0}^{M-1} b_k \rho_m^k
\end{equation*}
and
\begin{equation*}
\vv_m=
\left [
1 ,
\rho_m,
\rho_m^2,
\cdots,
\rho_{m}^{M-1}
\right ]^T~~(m=0,1,2,\dots,M-1),
\end{equation*}
where
\begin{equation*}
\rho_{m}=\text{e}^{{~\text{j}}2\pi m/M}~~(m=0,1,2,\dots,M-1)
\end{equation*}
are the $M$th complex roots of 1, and $\text{j}$ is the imaginary
unit. It can easily be verified that all eigenvectors are orthogonal
to each other, i.e.,
\begin{equation*}
\vv_{m_1} \perp  \vv_{m_2}~~\text{if}~~ m_1\neq m_2.
\end{equation*}

Matrix $B$ becomes a block circulant matrix when each $b_i$ is an $N\times N$ square matrix.
It has long been discovered that the computation of the eigenvalue and eigenvectors
of a block circulant matrix can also be simplified
by exploiting its cyclic symmetric
nature\cite{Friedman1961Eigenvalues}.
It can be shown that all the $M\times N$ eigenvector of the block circulant matrix
can be divided into $M$ groups, with each group containing $N$
vectors. The $n$th eigenvector in the $m$th group,
denoted as $\w_{m,n}$, should be of the following form
\begin{equation*}
\w_{m,n} =
\left[
1,
\rho_m ,
\rho_m^2 ,
\cdots,
\rho_m^{M-1}
\right]^T\vv_{m,n}
\end{equation*}
where $\vv_{m,n}$ is a vector of length $N$.
In order for $\w_{m,n}$ to be an eigenvector, the
following equation has to be satisfied
\begin{equation*}
B \w_{m,n} = \lambda_{m,n} \w_{m,n}
\end{equation*}
which can be expanded as follows
\begin{equation*}
\begin{array}{rcl}
(b_0 + \rho_m b_1+ \cdots + \rho_m ^{M-1}b_{M-1}) \vv_{m,n}& = & \lambda_{m,n} \vv_{m,n}\\
(b_{M-1} + \rho_m b_0+ \cdots + \rho_m ^{M-1}b_{M-2}) \vv_{m,n}& = & \lambda_{m,n}\rho_m  \vv_{m,n}\\
~&\vdots &~\\
(b_1 + \rho_m b_2+ \cdots + \rho_m ^{M-1}b_{0}) \vv_{m,n}& = & \lambda_{m,n}\rho_m^{M-1} \vv_{m,n}.
\end{array}
\end{equation*}
It can be verified that any two equations in the linear system above are
equivalent. Therefore $\lambda_{m,n}$ and $\vv_{m,n}$ can be found by
solving the first equation only. That is, they are
the $n$th eigenvalue and eigenvector of the $N\times N$ matrix $B_m$
\begin{equation*}
B_m=b_0 + \rho_m b_1+ \rho_m^2 b_2 + \cdots + \rho_m ^{M-1}b_{M-1}.
\end{equation*}
Therefore, the eigen analysis of the large matrix $B$ of dimension
$MN \times MN$ can be obtained via the eigen analysis of
$M$ smaller matrix $B_m$ of dimension $N\times N$, significantly
simplifying the computation.

\subsection{Jacobian matrix of the linearized RANS equations}
For a typical finite volume solver based on the RANS equations,
the vector $\U$ that satisfies
\begin{equation*}
\R(\U)=0
\end{equation*}
is the steady state solution. The time-dependent behavior of the flow
solution is governed by the unsteady form
\begin{equation*}
\dfrac{d\U}{dt}=\R(\U).
\end{equation*}
To study its linear stability, one could linearize the nonlinear residual vector
$\R$ about the steady state with respect to the flow solution vector $\U$ to
obtain the Jacobian matrix $A$
\begin{equation*}
A:=\df{\R}{\U}
\end{equation*}
and compute its eigenvalues and eigenvectors.

\begin{figure}[htb]
	\centering
	\includegraphics[width=0.55\textwidth]{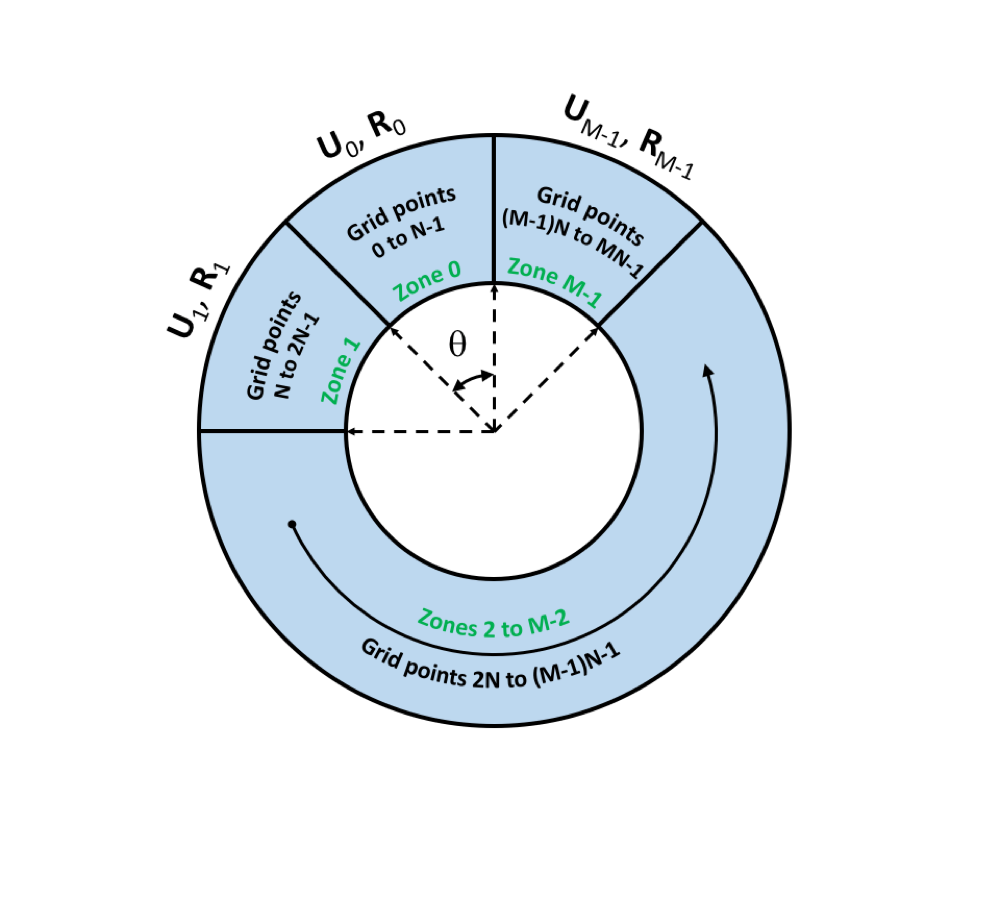}
	\caption{Arrangement of the computational meshes in the cyclic symmetric circular domain.}
	\label{circular-domain}
\end{figure}

To facilitate the discussion, it is assumed that the computational mesh
is cyclic symmetric with $M$-periodicity. Consequently, the
full-annulus mesh can be divided into $M$ non-overlapping zones, each of
which containing $N$ grid points, as illustrated in Fig.~\ref{circular-domain}.
In addition, the grid points in each zone is arranged in such a way that the
$n$th point in zone $m_1$ is the $n$th point in zone $m_2$ rotated by
an angle of $2(m_2-m_1)\pi /M$.

Due to the cyclic symmetry of the mesh, the Jacobian
matrix $A$ can be partitioned as a block matrix as follows
\begin{equation*}
A=
\left[
\begin{array}{cccc}
\df{\R_0}{\U_0} &\df{\R_0}{\U_1}  &\dots & \df{\R_0}{\U_{M-1}}\\
\df{\R_1}{\U_0} &\df{\R_1}{\U_1}  &\dots & \df{\R_1}{\U_{M-1}}\\
\vdots &\vdots&\ddots & \vdots\\
\df{\R_{M-1}}{\U_0} &\df{\R_{M-1}}{\U_1} &\dots & \df{\R_{M-1}}{\U_{M-1}}
\end{array}
\right]
\end{equation*}
where $\U_{m}$ and $\R_{m}$ are the vectors of the flow variables
and the nonlinear residual for the grid points in the $m$th zone, each
being a vector of length $5N$ for Euler/laminar or $6N$ for turbulent
flow with a one-equation turbulence model, and each block
$\df{\R_{[\cdot]}}{\U_{[\cdot]}}$ is a $5N\times 5N$ or $6N\times 6N$ matrix.

\begin{figure}[htb]
	\centering
	\includegraphics[width=0.55\textwidth]{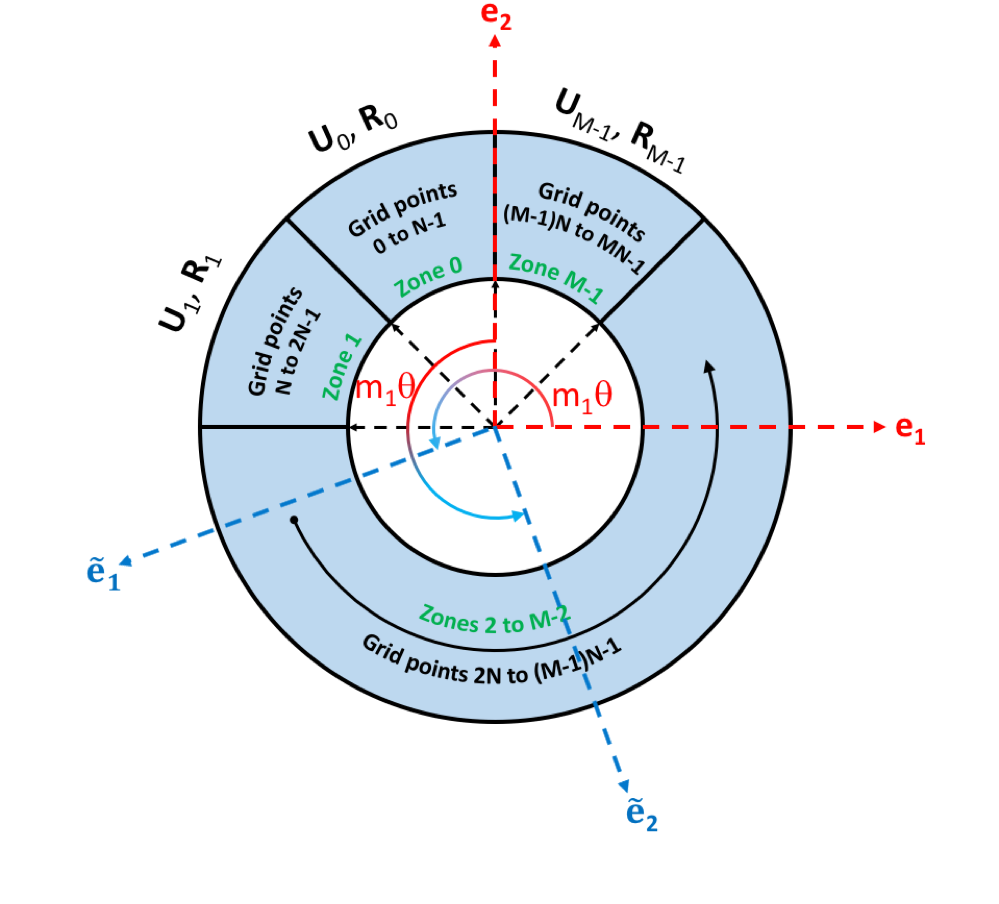}
	\caption{Arrangement of the computational meshes in the cyclic symmetric circular domain.}
	\label{rotation}
\end{figure}

Note that block matrices $\df{\R_{[\cdot]}}{ \U_{[\cdot]}}$ are
second-order tensors that are invariant to the coordinate system used, and
can be expressed using tensor product as
\begin{equation}
\df{\R_{m_1}}{ \U_{m_2}}
= \left[  M_{m_1,m_2} \right] _{i,j}\e_i \e_j
= \left[ \hat M_{m_1,m_2} \right] _{i,j}\hat \e_i \hat \e_j
\label{tensor}
\end{equation}
where $\{\e_1,\e_2,\e_3\}$ are the basis vectors in the coordinate frame fixed
to the computational mesh and $\e_i \e_j$ is short for the tensor product
$\e_i \otimes \e_j$. $M_{m_1,m_2}$ and $\hat M_{m_1,m_2}$
are the matrix representation of the second-order tensor $\df{\R_{m_1}}{\U_{m_2}}$
in the coordinate systems $\{\e_1,\e_2,\e_3\}$ and $\{\hat \e_1,\hat \e_2,\hat \e_3\}$
respectively. Let $\{\hat \e_1,\hat \e_2,\hat \e_3\}$ be a coordinate system
that is obtained by rotating $\{\e_1,\e_2,\e_3\}$ by $m_1\theta=2m_1\pi/M$, as shown in Fig.~\ref{rotation}
The basis vectors of the two coordinate systems can linked via
the transformation matrix $T^{m_1}$ via
\begin{equation*}
\e_i [T^{m_1}]_{i,j}= \hat \e_j,
\end{equation*}
where the matrix $T$ denotes the rotational transform of $\theta=2\pi/M$.
Substitute the $\hat \e_i$ in Eq.~\eqref{tensor} yields
\begin{equation*}
M_{m_1,m_2}=T^{m_1} \hat M_{m_1,m_2} T^{-m_1} .
\end{equation*}

In addition, it is straightforward to see that $\df{\R_0}{\U_{m_2-m_1}}$ viewed in $\{\e_1,\e_2,\e_3\}$
is identical to $\df{\R_{m_1}}{\U_{m_2}}$ viewed in $\{\tilde \e_1,\tilde \e_2,\tilde \e_3\}$. This means
the matrix representation of the former in $\{\e_1,\e_2,\e_3\}$ is that of the latter in
$\{\tilde \e_1,\tilde \e_2,\tilde \e_3\}$, i.e.,
\begin{equation*}
\hat M_{m_1,m_2}= M_{0,m_2-m_1}
\end{equation*}
and hence
\begin{equation*}
M_{m_1,m_2}= T^{m_1}M_{0,m_2-m_1}T^{-m_1}.
\end{equation*}
or equivalently,
\begin{equation*}
\df{\R_{m_1}}{\U_{m_2}}= T^{m_1}\df{\R_0}{\U_{m_2-m_1}}T^{-m_1}.
\end{equation*}
The Jacobian matrix $A$ becomes
\begin{equation*}
A=
\left[
\begin{array}{cccc}
\df{\R_0}{\U_0}        &   \df{\R_0}{\U_1}&\dots & \df{\R_0}{\U_{M-1}}\\
T \df{\R_0}{\U_{M-1}}     T^{-1} & T\df{\R_0}{\U_0}   T^{-1}&\dots &T \df{\R_0}{\U_{M-2}}  T^{-1}\\
\vdots &\vdots&\ddots & \vdots\\
T^{M-1} \df{\R_0}{\U_1} ˇ^{-(M-1)}&T^{M-1} \df{\R_0}{\U_{M-2}} T^{-(M-1)} &\dots & T^{M-1}\df{\R_0}{\U_0}T^{-(M-1)}
\end{array}
\right ]
\end{equation*}

Further more, change of variable is used for $\U_m$
\begin{equation*}
\tilde \U_m = T^{-m} \U_m
\end{equation*}
so that
\begin{equation*}
\df{\R_0}{\U_m}=\df{\R_0}{\tilde \U_m} T^{-m}
\end{equation*}

The Jacobian matrix $A$ can be further rewritten to be
\begin{align*}
A&=
\left[
\begin{array}{cccc}
  \df{\R_0}{\tilde \U_0}      &   \df{\R_0}{\tilde \U_1}T^{-1}&\dots &  \df{\R_0}{\tilde \U_{M-1}}T^{-(M-1)}\\
T\df{\R_0}{\tilde \U_{M-1}}& T\df{\R_0}{\tilde \U_0}T^{-1}&\dots &T  \df{\R_0}{\tilde \U_{M-2}}T^{-(M-1)}\\
\vdots &\vdots&\ddots & \vdots\\
T^{M-1}  \df{\R_0}{\tilde \U_1}  &T^{M-1}\df{\R_0}{\tilde \U_{M-2}}T^{-1} &\dots & T^{M-1} \df{\R_0}{\tilde \U_0}  T^{-(M-1)}
\end{array}
\right]
={ \mathcal {T}  } B { \mathcal {T}  }^{-1}\\
\end{align*}
with matrices $\mathcal{T}$ and $B$ defined as
\begin{equation*}
\mathcal{T}=	\left[
\begin{array}{cccc}
I & 0  & \cdots & 0\\
0 & T & \cdots  & 0\\
\vdots&\vdots &\dots & \vdots\\
0 &  0 & \cdots & T^{M-1}
\end{array}
\right],~~
B=\left[
\begin{array}{cccc}
\df{\R_0}{\tilde \U_0}                & \df{\R_0}{\tilde \U_1}        &\dots &  \df{\R_0}{\tilde \U_{M-1}}     \\
\df{\R_0}{\tilde \U_{M-1}}          & \df{\R_0}{\tilde \U_0}       &\dots &  \df{\R_0}{\tilde \U_{M-2}}     \\
\vdots &\vdots&\ddots & \vdots\\
  \df{\R_0}{\tilde \U_1}              &  \df{\R_0}{\tilde \U_{M-2}}&\dots & \df{\R_0}{\tilde \U_0}
\end{array}
\right]
\end{equation*}

Therefore, the Jacobian matrix $A$ is similar to a block-circulant matrix $B$
and thus they have the same eigenvalues. To compute the eigenvalues and
eigenvectors of the matrix $B$, we first compute $\lambda_{m,n}$ and
$\vv_{m,n}$, $n$th eigenvalue and eigenvector of the $N\times N$
matrix $B_m$
\begin{equation*}
B_{m}=
\df{\R_0}{\tilde \U_0}
+
\rho_m  \df{\R_0}{\tilde \U_1}
+
\rho_m^2  \df{\R_0}{\tilde \U_2}
+
\cdots
+
\rho_m ^{M-1} \df{\R_0}{\tilde \U_{M-1}}  ,
\end{equation*}
and the eigenvectors of $B$, $\w_{m,n}$, can be computed as
\begin{equation*}
\w_{m,n} =
\left[1,
\rho_m,
\rho_m^2 ,
\cdots ,
\rho_m^{M-1}\right]^T \vv_{m,n}   .
\end{equation*}
Finally, the eigenvectors of the Jacobian matrix $A$ are $\mathcal {T}\w_{m,n}$,
which can be calculated as follows
\begin{equation*}
\mathcal {T}\w_{m,n}
=
\left[
1,
\rho_m T,
\rho_m^2  T^2,
\cdots ,
\rho_m^{M-1}  T^{M-1}
\right]^T \vv_{m,n}   .
\end{equation*}

{ It can be seen that for any of the $M$ eigenvectors of $B_m$, $\vv_{m,n} (m=0,1,2,\dots, M-1)$,
the corresponding $M$ eigenvectors of $A$, $ \mathcal {T}\w_{m,n} (m=0,1,2,\dots, M-1)$
are modes with different nodal diameters. Same as for the circulant block matrix, to compute
the eigenvalues and eigenvectors of the Jacobian matrix $A$, one only needs to solve the
eigenvalue problem of dimension $N\times N$, significantly reducing the computational cost,
especially for large CFD analysis.

\section{Results}
\label{sec:result}
The test case used to demonstrate the usefulness of the proposed
method is an annulus compressor cascade case.
The cascade is produced by taking the surface of revolution at the
50\% height of the NASA Rotor~67~\cite{strazisar1989laser}. The
three-dimensional mesh has one cell in the radial direction and the
subsequent analysis is thus a quasi-3D one. The compressor
cascade rotates at 16,043 Revolutions per minute. The annular cascade
has 22 passages in total. The mesh for the whole annulus is produced
by meshing a single passage first and then copying the mesh to the
whole annulus. The mesh for the entire domain has 126,808 grid points,
with an average of 5,764 per zone.

\subsubsection{Steady state calculation}
The turbomachinery nonlinear flow solver NutsCFD is used to compute the
performance of the cascade. The speedline, computed
by gradually raising the back pressure until the steady state can no
long converge, is shown in Fig.~\ref{speedline}.
The fully-converged steady state flow solution corresponding to
the condition marked in red in Fig.~\ref{speedline} (to be called
condition `A' in the remainder of the paper) is visualized in
Fig.~\ref{flow}, where the shock/boundary layer interaction
can clearly be seen. All steady state calculations are performed
using a single passage. The visualized flow fields are produced
by copying the solution in the single passage to the whole annulus.

\begin{figure}[htb]
	\centering
	\includegraphics[width=0.9\textwidth]{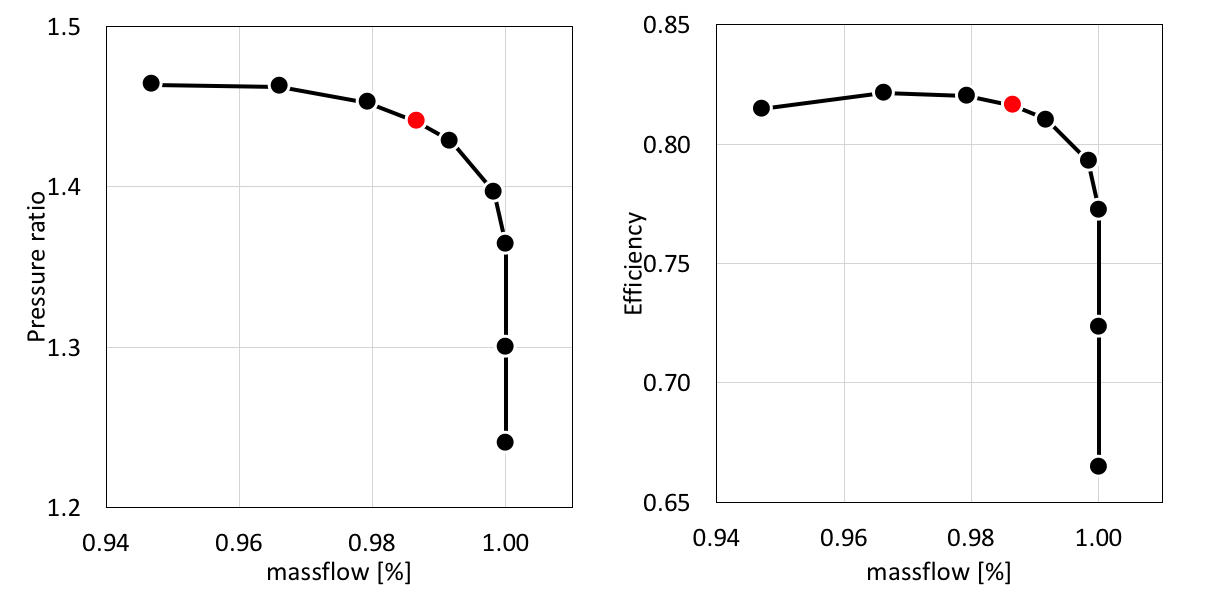}
	\caption{Total pressure ratio and isentropic efficiency of the cascade with the
	red circle denoting the condition for which eigenvalue analysis is performed.}
	\label{speedline}
\end{figure}

\begin{figure}[htb]
	\centering
	\includegraphics[width=0.45\textwidth]{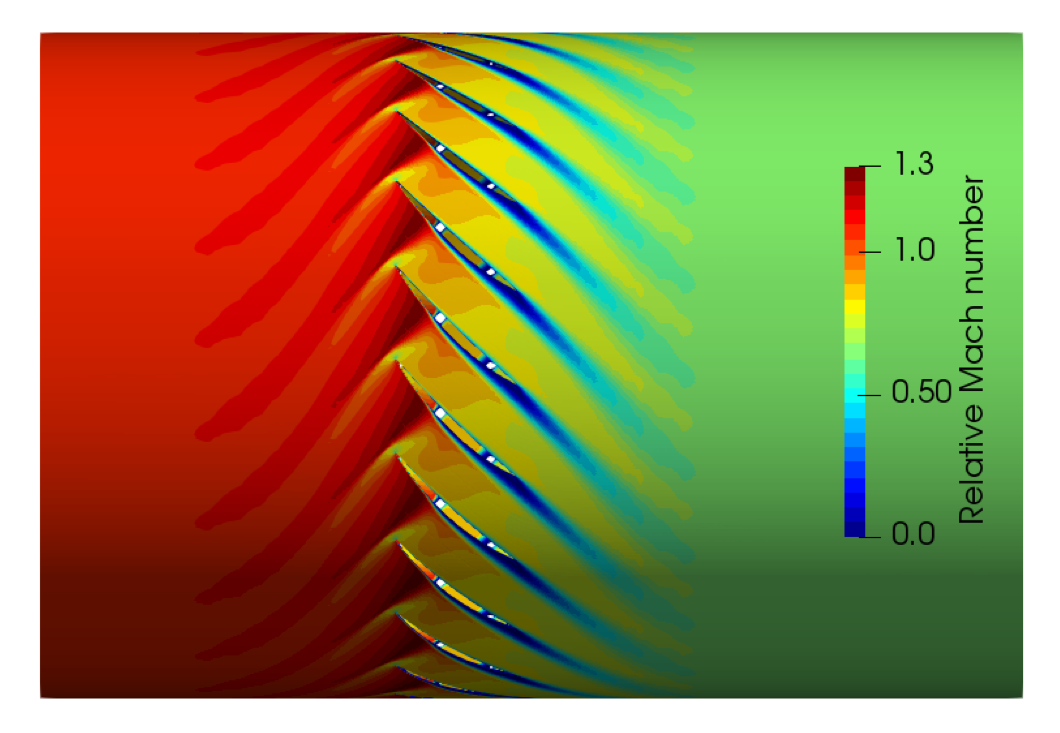}
	~~~~~
	\includegraphics[width=0.45\textwidth]{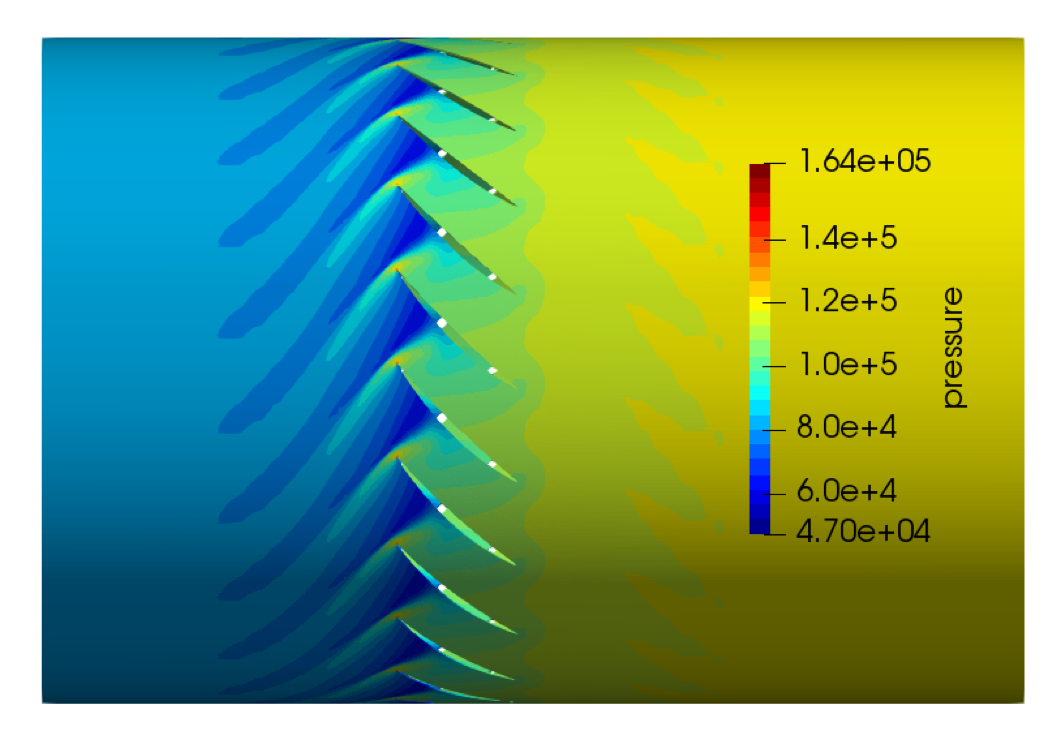}
	\caption{Contour plots for the relative Mach number (left)
		and pressure (right) for the converged steady state solution
		for condition `A'.}
	\label{flow}
\end{figure}

\subsubsection{Eigenvalue analysis}
To perform the eigenvalue analysis for condition `A', we first compute
the Jacobian matrix for the single passage $\df{\R_0}{\tilde \U_0}$
and then exploiting the periodic boundary condition in the steady
solver to obtain $\df{\R_0}{\tilde \U_1} $ and
$\df{\R_0}{\tilde \U_{M-1}} $. Note that for RANS flow solver,
$\df{\R_0}{\tilde \U_m},  (m=2,3,\dots, M-2) $ are all zeros, so
$B_m$ can be simplified as
\begin{equation}
B_{m}=\df{\R_0}{\tilde \U_0}+\rho_m  \df{\R_0}{\tilde \U_1} +\rho_m ^{M-1} \df{\R_0}{\tilde \U_{M-1}}.
\label{eq:Bm}
\end{equation}
It is well known both numerically and experimentally that the least stable modes
when instability appears in rotating flows, the characteristic frequency of the
instability is close to the rotating frequency, which is $\dfrac{2\pi \times 16,043}{60}=1680 Hz$
for the test case used. Consequently, the matrix $B_m$ is scaled by $1/1680$ so
that the imaginary part of the eigenvalue is equivalent to the frequency of the
corresponding eigenmode expressed in terms of engine order (EO).
Since the Jacobian matrix is of moderate size (of dimension $34,584\times 34,584$, and
with nearly $3.7$ million non-zero entries), the sparse matrix eigenvalue solver, ``eigs'',
in SciPy~\cite{jones2001scipy} is use on a single-processor computer. The ``eigs''
function in SciPy is a wrapper to the widely used eigenvalue
solver library ARPACK~\cite{lehoucq1998arpack}.

Two different methods (listed in Tab.~\ref{tab:methods}) are
used to compute the eigenvalues near the imaginary axis.
Method 1 computes the eigenvalues of the Jacobian matrix
for the whole annulus of 22 sectors, $A$. We apply complex shifts of
$\sigma=i, 2i$ and $3i$ and 30 interior eigenvalues are computed
for each shift. A total of 90 eigenvalues can computed and they are
shown in Fig.~\ref{eigenvalue} as red crosses.

Method 2 computes the eigenvalues using 1 sector only. Nevertheless, it admits
eigenvalues for the whole annular domain by setting $\rho_m$ in Eq.~\eqref{eq:Bm}
to $\rho_{m}=\text{e}^{{~\text{j}}2\pi m/M}~~(m=0,1,2,\dots,M-1)$. Same as method 1,
interior eigenvalues are computed using three different shifts for each $\rho_m$. For
each shift, 2 eigenvalues are computed, resulting in a total of 132 eigenvalues
(some of them appear more than once for different shifts). All eigenvalues
computed using method 2 are also plotted in Fig.~\ref{eigenvalue} as blue dots.
If can be seen that the least stable eigenvalues (closest to the imaginary axis) are
located around $(0,1)$ and a zoomed view of the eigenvalues in the nearly region
is shown on the right side in Fig.~\ref{eigenvalue}. It can be seen that both methods
can capture a relevant subset of the eigenvalues and are both sufficient in order to
study the linear stability of the underlying flow problem. However, method 2 is
based only on one sector of the whole annulus domain, and thus has a significantly
lower memory overhead and also costs much less CPU time than method 1.

Another advantage of method 1, besides the lower computational cost, is that
each eigenvalue/eigenvector computed has a known nodal diameter, which is
solely determined by the $m$ value it is computed with. The nodal diameter
information reveals the circumferential spatial periodicity of each particular
eigenmode.

\begin{table}[htb]
	\centering
	\caption{Two methods for computing the spectrum.}
	\label{tab:methods}
	\begin{tabular}{ c c c c }
		\hline \hline
		method & number of sectors & matrix dimension & $\theta$ value \\ \hline
        1 & 22&  $760,848 \times 760,848$& no need to specify \\
		2 & 1  &  $34,584  \times 34,584$ & ${2m\pi}/{M}~(m=0,1,2,\dots,M-1)$ \\
		 \hline \hline
	\end{tabular}
\end{table}

\begin{figure}[htb]
\centering
\includegraphics[width=0.45\textwidth]{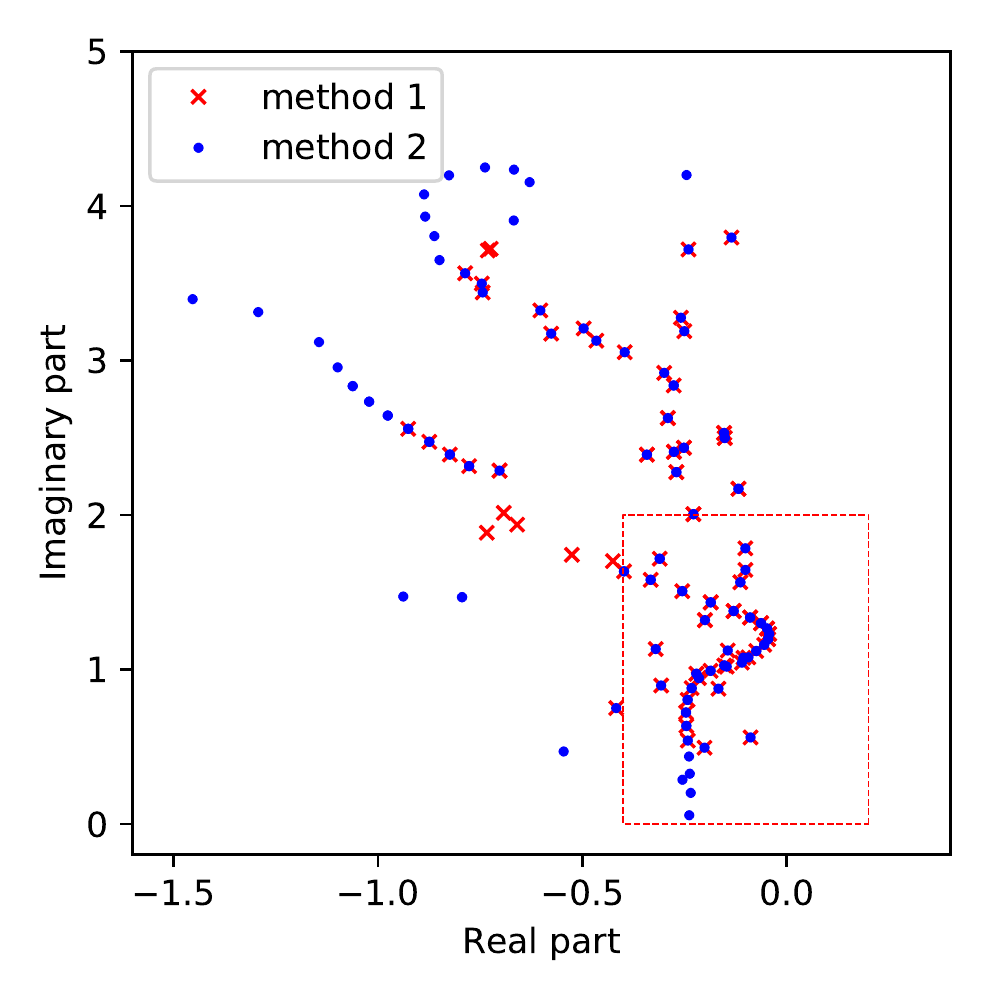}
\includegraphics[width=0.45\textwidth]{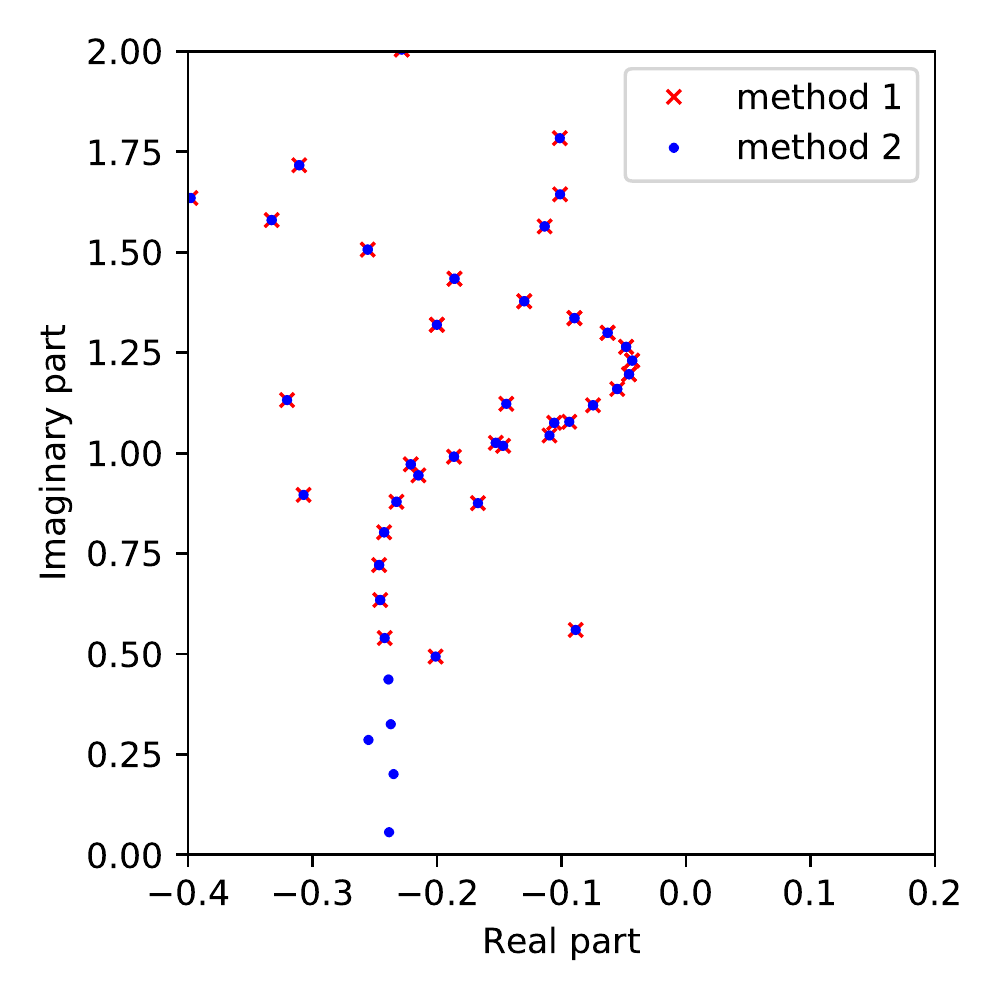}
\caption{Computed spectrum near the imaginary axis for conditions `A'.
Red cross: eigenvalues computed with method 1 using 22 sectors;
blue dot: eigenvalues computed with method 2 using 1 sector only.
A zoomed view of the eigenvalues in the region in the red box is
shown on the right.}
\label{eigenvalue}
\end{figure}

\section{Conclusions}
\label{sec:conclusion}
A method to simplify the eigenvalue analysis for flow problems in
circular cyclic symmetric domains are proposed. The method first
transform the Jacobian matrix for the whole annulus into a block
circulant matrices, and subsequently reduces the eigenvalue
problem of the whole circular domain to one sector only. This
reduces the size of the original eigenvalue problem by a factor of
$M$, $M$ being the number of passages/blades. The proposed
method is applied to the eigenvalue analysis for an annular
compressor cascade with $22$ sectors, and it is shown that by
using $1$ sector only, spectrum of the whole domain can be
obtained at a much lower memory as well as CPU time cost,
demonstrating the advantage of the proposed method. Future work
will explore the application of the method to realistic
three-dimensional compressors to study the flow instability problems.
In addition, the eigenvectors obtained can be used to
construct reduced-order models for parametric study and optimization
for rotating-flow instability.

\bibliography{xu}

\end{document}